\documentclass[10pt]{amsart}

\usepackage{Settings}

\begin{document}

\title{A note on Duality and the Atiyah--Hirzebruch spectral sequence}
\date{March 19, 2026}

\author{Maximilian David Hans}
\address{School of Mathematical Sciences, University of Southampton, United Kingdom}
\email{maximilian.hans@soton.ac.uk}

\begin{abstract}
We show that, for a finite spectrum $X$, Spanier--Whitehead duality induces an isomorphism between the cohomological and homological Atiyah--Hirzebruch spectral sequences. As an application, it follows that Poincar\'e duality for a Poincar\'e duality complex that is oriented over a ring spectrum $\cR$ induces an isomorphism between the two spectral sequences.
\end{abstract}

\maketitle


\section{Introduction}

Spanier--Whitehead duality assembles into a duality functor $\bD\colon (\Sp^\fin)^\op\ra \Sp^\fin$ on finite spectra. For an arbitrary spectrum $\cE$, and a finite spectrum $X$, the identification $\bD(X)\simeq \F(X,\S)$ yields the Spanier--Whitehead duality isomorphism $\SW\colon \cE^k(X)\ra \cE_{-k}(\bD(X))$ by adjunction. A CW filtration on $X$ induces a dual CW filtration on $\bD(X)$. We investigate the interplay between the Spanier--Whitehead duality isomorphism and the Atiyah--Hirzebruch spectral sequence using these filtrations. The result may be known to experts, but we were unable to locate a precise reference in the literature.

\begin{MainThm}\label{thm:A}
    Let $\cE$ be a spectrum, and $X$ a finite spectrum. The Spanier--Whitehead duality isomorphism $\SW\colon\cE^*(X)\ra \cE_{-*}(\bD(X))$ induces an isomorphism between the cohomological and homological Atiyah--Hirzebruch spectral sequences computing $\cE^*(X)$ and $\cE_{-*}(\bD(X))$ respectively. More precisely, there is an isomorphism $\smash{(E^{p,q}_r,d_r)\ra (E_{-p,-q}^r,d^r)}$ on every page of the two spectral sequences induced by Spanier--Whitehead duality.
\end{MainThm}

The unital and multiplicative structure provided by a ring spectrum $\cR$ naturally induces a notion of $\cR$-orientability for a $d$-dimensional Poincar\'e duality complex $X$. Briefly speaking, an $\cR$-orientation is a globally consistent way to assign to each fibre of the Spivak normal fibration of $X$ a specified unit in $\pi_0(\cR)$. An $\cR$-orientation in turn yields a Poincar\'e duality isomorphism $\smash{\PD\colon \cR^k(X)\ra \cR_{d-k}(X)}$ between cohomology and homology. Using the fact that Poincar\'e duality factors as the composition of the Thom isomorphism followed by the Spanier--Whitehead duality isomorphism, we use \cref{thm:A} to show that the Poincar\'e duality isomorphism extends to an isomorphism between the cohomological and homological Atiyah--Hirzebruch spectral sequences. This allows us to leverage knowledge of differentials in the cohomological spectral sequence to obtain differentials in the homological version. For manifolds, this result was previously shown by Vick \cite{VickPDAHSS} using properties of the slant product, and in the case of $\MSpin$ and spin manifolds by Davis \cite[{Lemma~3.8}]{DavisBorelNovikov}.

\begin{MainThm}\label{thm:B}
    Let $\cR$ be a ring spectrum, and $X$ a $d$-dimensional Poincar\'e duality complex oriented over $\cR$. Poincar\'e duality for the ring spectrum $\cR$ induces an isomorphism between the cohomological and homological Atiyah--Hirzebruch spectral sequences computing $\cR^*(X)$ and $\cR_{d-*}(X)$ respectively. More precisely, there is an isomorphism $\smash{(E^{p,q}_r,d_r)\ra (E_{d-p,-q}^r,d^r)}$ on every page of the two spectral sequences induced by Poincar\'e duality.
\end{MainThm}

\section{Preliminaries and conventions}\label{sec:preliminaries}

We begin by briefly introducing the notation and standard results that will be used throughout this article. This includes a discussion on spectra, which is important for \cref{sec:SW-spectra} and \cref{sec:SW-AHSS}, while standard results on Poincar\'e duality complexes are used only in \cref{sec:PD-Orientation}.

\subsection*{Spectra}

Throughout this article, we work with the $\infty$-category of spaces $\sS$, and the stable $\infty$-category of spectra $\Sp$, defined as the limit
\[
    \Sp\coloneqq\lim(\dots\xrightarrow{\Omega}\sS_*\xrightarrow{\Omega}\sS_*)
\]
taken in the $\infty$-category of $\infty$-categories (\cite[{Section~1.4.2,~Definition~1.4.3.1}]{HigherAlgebra}).
This may be safely ignored by readers who are only familiar with spectra in the classical context, and mostly serves as a framework of consistent notation. This subsection provides a summary of standard definitions and results, presented in the language and framework established in \cite{HigherAlgebra}.

\smallskip
\noindent\textit{Shifting spectra.}
The adjoint functors
\[\begin{tikzcd}
	{\Sigma\colon\Sp} & {\Sp\flippedcolon\Omega}
	\arrow[shift left, from=1-1, to=1-2]
	\arrow[shift left, from=1-2, to=1-1]
\end{tikzcd}\]
are mutual inverses, and are equivalent to shifting the sequence of spaces given by a spectrum. For a spectrum $\cE$, we write $\Sigma^n \cE \eqqcolon \cE[n]$ and $\Omega^n\cE\eqqcolon \cE[-n]$.

\smallskip
\noindent\textit{Suspension spectra.}
The functor $\Sigma^\infty\colon \sS_*\ra \Sp$, which on objects maps a space to its suspension spectrum, fits into an adjunction
\[\begin{tikzcd}
	{\Sigma^\infty\colon\sS_*} & {\Sp\flippedcolon\Omega^\infty}
	\arrow[shift left, from=1-1, to=1-2]
	\arrow[shift left, from=1-2, to=1-1]
\end{tikzcd}\]
between the $\infty$-category of pointed spaces and spectra. For a pointed space $X$, we generally write $\Sigma^\infty X$ for its suspension spectrum. If $X$ is not pointed, we consider $X_+$, the space $X$ together with a disjoint base point before passing to the suspension spectrum.

\smallskip
\noindent\textit{(Co)homology groups valued in spectra.}
\begin{enumerate}
    \item[(1)] For a spectrum $\cE$, and a spectrum $X$, the \emph{$\cE$-homology groups} of $X$ are defined as $\cE_n(X)\coloneqq \pi_n(\cE\otimes X)$, where the \emph{smash product of spectra} $-\otimes-\colon\Sp\times \Sp\ra \Sp$ is uniquely determined by the fact that it preserves small colimits in both variables, and that the unit object of $\Sp$ is the sphere spectrum. In particular, this endows $\Sp$ with a symmetric monoidal structure, which is unique up to contractible choice (\cite[{Corollary~4.8.2.19}]{HigherAlgebra}). In the case where $X$ is an unpointed space, we write $\cE_n(X)\coloneqq \pi_n(\cE\otimes \Sigma^\infty X_+)$, and in the case of $X$ being pointed, define the \emph{reduced} $\cE$-homology groups by $\wt{\cE}_n(X)\coloneqq \pi_n(\cE\otimes \Sigma^\infty X)$. In particular, for an unpointed space $X$, there is an equality $\cE_n(X)=\wt{\cE}_n(X_+)$.
    \item[(2)] The \emph{$\cE$-cohomology groups} of a spectrum $X$ are defined as $\cE^n(X)\coloneqq \pi_{-n}(\F(X,\cE))$, where the spectrum $\F(\Sigma^\infty X_+,\cE)$ is the \emph{mapping spectrum} from $\Sigma^\infty X$ to $\cE$. Its associated infinite loop space is the mapping space $\Map_{\Sp}(X,\cE)$. This leads to the classical adjunction $\F(A,\F(B,C))\simeq \F(A\otimes B, C)$ for spectra $A$, $B$, and $C$, since the smash product functor is left adjoint to the internal $\Hom$ functor. In the case where $X$ is an unpointed space, we write $\cE^n(X)\coloneqq \pi_{-n}(\F(\Sigma^\infty X_+,\cE))$, and in case $X$ is already pointed, the \emph{reduced} $\cE$-cohomology groups are defined as $\wt{\cE}^n(X)\coloneqq \pi_{-n}(\F(\Sigma^\infty X,\cE))$.
\end{enumerate}
The case of $\cE$-(co)homology of the sphere spectrum $\S\coloneqq \Sigma^\infty\S^0$ is particularly simple---we obtain the homotopy groups of $\cE$. We abbreviate $\cE_n(\S)\simeq \cE_n(*)\eqqcolon \cE_n$ and $\cE^n(\S)\simeq \cE^n(*)\eqqcolon \cE^{n}$, where we assume the space $*$ to be unpointed. In particular, there is an isomorphism of abelian groups between $\cE^n$ and $\cE_{-n}$.

\smallskip
\noindent\textit{Homotopy classes of maps.}
The homotopy classes of maps between two spectra $A$ and $B$ are defined as $[A,B]\coloneqq \pi_0(\Map_{\Sp}(A,B))$, where the mapping space $\Map_{\Sp}(A,B)$ of spectra is an object in $\sS$. This is equivalent to the group $\pi_0(\F(A,B))$ via the chain of equivalences
\[
    \pi_0(\F(A,B))\simeq \pi_0(\Map_{\Sp}(\S,\F(A,B)))\simeq \pi_0(\Map_{\Sp}(\S\otimes A,B))\simeq\pi_0(\Map_{\Sp}(A,B)).
\]

\subsection*{Poincar\'e duality complexes}

A $d$-dimensional Poincar\'e duality complex is a connected space $X$ which has the homotopy type of a finite CW complex, such that there exists a pair $(\l,[X])$ where $\l$ is the orientation local system (locally infinite cyclic), and $[X]$ is a class in $\H_n(X;\l)$ such that the map
\[
    \H^k(X;\m)\xlongrightarrow{-\cap[X]}\H_{d-k}(X;\m\otimes \l)
\]
is an isomorphism for any local system of abelian groups $\m$ on $X$ for all $k\in \Z$. We call the class $[X]$ the \emph{fundamental class of $X$}.
Poincar\'e duality complexes have been thoroughly studied and play a prominent r\^{o}le in surgery theory. We refer the reader to \cite[{Section~9}]{RanickiSurgery} for basic properties of (geometric) Poincar\'e duality complexes, and \cite[{Appendix~A}]{ReducibilityLand} for a view on Poincar\'e duality complexes via parametrised homotopy theory.

\smallskip
\noindent\textit{Spivak normal fibration.} 
Let $\BG(r)$ be the classifying space of homotopy equivalences of the $(r-1)$-sphere, and let $\BG\coloneqq \colim_{r}\BG(r)$. Note that $\BG(r)$ classifies spherical fibrations of rank $r$ (the fibre being an $(r-1)$-sphere), whereas $\BG$ classifies stable spherical fibrations.
Consider a finite CW complex $X$ of dimension $d$ with an embedding $i\colon X\hookrightarrow \R^\ell$ into high-dimensional Euclidean space. A thickening $N(X)$ of the image is homotopy equivalent to $X$, and thus we may view the boundary inclusion $\partial N(X)\hookrightarrow N(X)$ as a map $\partial N(X)\ra X$. Spivak showed in \cite{Spivak} that the fibre of this map has the homotopy type of a sphere $\S^{\ell-d-1}$ if and only if $X$ is a Poincar\'e duality complex. In this case, stabilisation via embeddings $\R^\ell\hookrightarrow \R^{\ell+1}\hookrightarrow\R^{\ell+2}\hookrightarrow\dots$ yields a stable spherical fibration $\SF(X)\colon X\ra\BG$ called the \emph{Spivak normal fibration} of $X$. For Poincar\'e duality complexes, the Spivak normal fibration plays the role of the stable normal bundle for manifolds.

\smallskip
\noindent\textit{Thom spaces.} 
 For a spherical fibration $\zeta\colon X\ra \BG(r)$ of rank $r$, there is a well-defined Thom space via $\Th(\zeta)\coloneqq \cofib(S(\zeta)\ra X)$, where $S(\zeta)$ is the total space of the spherical fibration. Note that $\Th(\zeta)$ is naturally pointed. In particular, if $J\colon \BO(r)\ra \BG(r)$ is the $J$-homomorphism, which sends a vector bundle to its underlying spherical fibration by restriction, and $\xi\colon X\ra \BO(r)$ is a vector bundle of rank $r$, there is a homotopy equivalence $\Th(\xi)\simeq \Th(J\circ \xi)$ (see for example \cite[{Proposition~9.24}]{RanickiSurgery}). Hence, the homotopy type of the Thom space of a vector bundle only depends on the underlying spherical fibration.

\section{Spanier--Whitehead duality for spectra}\label{sec:SW-spectra}

This section serves as a recollection of well-known facts about Spanier--Whitehead duality that we require in the proof of \cref{thm:A}. Furthermore, we provide the necessary setup of the Atiyah--Hirzebruch spectral sequence for said proof, which we tend to in \cref{sec:SW-AHSS}.

\subsection*{Spanier--Whitehead duality and (co)fibrations}

A spectrum $X$ is dualisable if there exists a spectrum $\bD(X)$ such that $\bD(X)\otimes -$ is right adjoint to $X\otimes-$. This occurs if and only if $X$ is a finite spectrum, in which case its \emph{Spanier--Whitehead dual} is given by $\bD(X)\simeq \F(X,\S)$. We may take this as the definition of the Spanier--Whitehead dual of a finite spectrum. Restricting to the full $\infty$-subcategory $\smash{\Sp^\fin}$ of finite spectra, we obtain a duality functor $\bD\colon \smash{(\Sp^\fin)^\op}\ra \smash{\Sp^\fin}$. 
It is well-known that, since the $\infty$-category of spectra is stable, this functor preserves the zero object as well as fibre sequences as shown in \cref{lem:cofib}---in other words, the Spanier--Whitehead duality functor $\bD\colon \smash{(\Sp^\fin)^\op}\ra \smash{\Sp^\fin}$ is exact.

The interaction between a finite spectrum $X$ and its Spanier--Whitehead dual $\bD(X)$ is particularly simple through the lens of (co)homology, which is the driving factor behind \cref{thm:A}. For an arbitrary spectrum $\cE$, adjunction yields the equivalence
\begin{align*}
         \cE^{k}(X) & \simeq[\S[-k]\otimes X, \cE] \\
         & \simeq [\S[-k], \bD(X)\otimes \cE] \\
         & =\cE_{-k}(\bD(X))
\end{align*}
and we denote the isomorphism by $\SW\colon \cE^k(X)\ra \cE_{-k}(\bD(X))$. This natural isomorphism is called the \emph{Spanier--Whitehead duality isomorphism}, evidently induced by the functor $\bD\colon \smash{(\Sp^\fin)^\op}\ra \smash{\Sp^\fin}$ and its defining property. Exactness of this functor is crucial since the exact couple defining the Atiyah--Hirzebruch spectral sequence is built from cofibrations. Therefore, we include an elaborate discussion below. Note that both $\Sp$ and $\Sp^\fin$ are stable $\infty$-categories, so we may freely interchange the words \emph{fibration} and \emph{cofibration}.

\begin{Lem}\label{lem:cofib}
    Let $f\colon X\ra Y$ be a map of finite spectra, together with an associated (co)fibre sequence
\[
    X\xlongrightarrow{f}Y\lra\cofib(f).
\]
    The Spanier--Whitehead dual of $\cofib(f)$ satisfies the natural equivalence $\bD(\cofib(f))\simeq \fib(\bD f)$.
\end{Lem}

\begin{proof}
    Since $\bD(-)\coloneqq \F(-,\S)$, the composition $\bD(\cofib(f))\ra \bD(Y)\ra \bD(X)$ comes with a null-homotopy induced by the data of the fibre sequence $X\ra Y\ra \cofib(f)$. The induced maps are obtained by pre-composition. By the universal property of limits, we therefore obtain a map $\bD(\cofib(f))\ra \fib(\bD f)$. The data is collected in the diagram
\[\begin{tikzcd}
	{\bD(\cofib(f))} \\
	{\fib(\bD f)} & {\bD(Y)} & {\bD(X)}
	\arrow[dashed, from=1-1, to=2-1]
	\arrow[from=1-1, to=2-2]
	\arrow[from=2-1, to=2-2]
	\arrow["\bD f", from=2-2, to=2-3]
\end{tikzcd}\]
    where the lower sequence is the fibre sequence associated to the map $\bD f$. Let $T$ be an arbitrary test spectrum. The vertical equivalences in the diagram
\[\begin{tikzcd}
	{[T\otimes Y[1],\S]} & {[T\otimes X[1],\S]} & {[T\otimes \cofib(f),\S]} & {[T\otimes Y,\S]} & {[T\otimes X,\S]} \\
	{[T,\bD(Y)[-1]]} & {[T,\bD(X)[-1]]} & {[T,\bD(\cofib(f))]} & {[T,\bD(Y)]} & {[T,\bD(X)]}
	\arrow[from=1-1, to=1-2]
	\arrow["\simeq", from=1-1, to=2-1]
	\arrow[from=1-2, to=1-3]
	\arrow["\simeq", from=1-2, to=2-2]
	\arrow[from=1-3, to=1-4]
	\arrow["\simeq", from=1-3, to=2-3]
	\arrow[from=1-4, to=1-5]
	\arrow["\simeq", from=1-4, to=2-4]
	\arrow["\simeq", from=1-5, to=2-5]
	\arrow[from=2-1, to=2-2]
	\arrow[from=2-2, to=2-3]
	\arrow[from=2-3, to=2-4]
	\arrow[from=2-4, to=2-5]
\end{tikzcd}\]
    are induced by adjunction. In particular, the left-most equivalence, which switches the sign of the shift, stems from the fact that $\F(-,\S)$ takes pushout diagrams to pullback diagrams. From the previous discussion, we obtain a comparison of long exact sequences
\[\begin{tikzcd}
	{[T,\bD(Y)[-1]]} & {[T,\bD(X)[-1]]} & {[T,\bD(\cofib(f))]} & {[T,\bD(Y)]} & {[T,\bD(X)]} \\
	{[T,\bD(Y)[-1]]} & {[T,\bD(X)[-1]]} & {[T,\fib(\bD f)]} & {[T,\bD(Y)]} & {[T,\bD(X)]}
	\arrow[from=1-1, to=1-2]
	\arrow[equals, from=1-1, to=2-1]
	\arrow[from=1-2, to=1-3]
	\arrow[equals, from=1-2, to=2-2]
	\arrow[from=1-3, to=1-4]
	\arrow[dashed, from=1-3, to=2-3]
	\arrow[from=1-4, to=1-5]
	\arrow[equals, from=1-4, to=2-4]
	\arrow[equals, from=1-5, to=2-5]
	\arrow[from=2-1, to=2-2]
	\arrow[from=2-2, to=2-3]
	\arrow[from=2-3, to=2-4]
	\arrow[from=2-4, to=2-5]
\end{tikzcd}\]
    and the central map is an isomorphism by the Five Lemma. The ($\infty$-categorical) Yoneda Lemma yields the desired equivalence $\bD(\cofib(f))\simeq \fib(\bD f)$. Since it is induced by adjunction, this equivalence is natural.
\end{proof}

\subsection*{CW filtrations on Spanier--Whitehead duals}

By definition, a finite spectrum $X$ admits a CW filtration  $0= X^{(-\infty)}\ra\dots\ra 0 \ra X^{(m)}\ra\dots \ra X^{(n)}\simeq X$ of finite spectra for specific integral bounds $m$ and $n$. The given maps $i\colon X^{(p-1)}\ra X^{(p)}$ are characterised by the fact that the cofibre admits an equivalence $\cofib(i)\simeq \smash{\bigoplus_{I_p}\S[p]}$ to a sum of shifted sphere spectra. Note that $X$ is necessarily $(m-1)$-connective. Suppressing the bounds, we write $\smash{X^{(*)}}$ for a chosen CW filtration on a finite spectrum $X$. The Spanier--Whitehead duality functor $\bD\colon \smash{(\Sp^\fin)^\op}\ra \smash{\Sp^\fin}$ does not preserve a chosen CW filtration, but we obtain a CW filtration on $\bD(X)$ as follows.

\begin{Def}\label{def:dual-cw}
    Let $X$ be a finite spectrum together with a chosen CW filtration $X^{(*)}$. The \emph{dual CW filtration} is a CW filtration on the Spanier--Whitehead dual $\bD(X)$ by setting $\bD(X)^{(-p)}\coloneqq \fib(\bD(X)\ra \bD(X^{(p-1)}))$.
\end{Def}

\begin{Lem}\label{lem:dual-cw}
    Let $X$ be a finite spectrum together with a chosen CW filtration $X^{(*)}$ with integral bounds $m$ and $n$ as above. The induced dual CW filtration as defined in \cref{def:dual-cw} is a well-defined finite CW filtration on $\bD(X)$.
\end{Lem}

\begin{proof}
    Let $X$ be a finite spectrum with a chosen CW filtration, and consider the fibre sequence
\[
        X^{(p-1)}\xlongrightarrow{i}X^{(p)}\xlongrightarrow{j}\cofib(i)\eqqcolon X^{(p)}/X^{(p-1)}.
\]
    By \cref{lem:cofib}, the Spanier--Whitehead dual of $X^{(p)}/X^{(p-1)}$ admits a natural equivalence to the fibre $\fib(\bD i)\simeq\cofib(\bD i)[-1]\eqqcolon\bD(X^{(p-1)})/\bD(X^{(p)})[-1]$. Since the inclusion $i\colon X^{(p-1)}\ra X^{(p)}$ comes from the CW filtration on $X$, we have an equivalence $X^{(p)}/X^{(p-1)}\simeq \bigoplus_{I_p}\S[p]$, and thus
\[
    \bD(X^{(p-1)})/\bD(X^{(p)})[-1]\simeq \bD(X^{(p)}/X^{(p-1)})\simeq \bD(\sideset{}{_{I_p}}\bigoplus\S[p])\simeq \F(\sideset{}{_{I_p}}\bigoplus\S[p],\S)\simeq \sideset{}{_{I_p}}\bigoplus\S[-p].
\]
    By definition of the dual CW filtration, we have a commutative square
\[\begin{tikzcd}[cramped]
	{\bD(X)^{(-p-1)}} & {\bD(X)} & {\bD(X^{(p)})} \\
	{\bD(X)^{(-p)}} & {\bD(X)} & {\bD(X^{(p-1)})} \\
	{\bD(X)^{(-p)}/\bD(X)^{(-p-1)}} & 0 & {\bD(X^{(p-1)})/\bD(X^{(p)})}
	\arrow[from=1-1, to=1-2]
	\arrow[from=1-1, to=2-1]
	\arrow[from=1-2, to=1-3]
	\arrow["\id", from=1-2, to=2-2]
	\arrow["{\bD i}", from=1-3, to=2-3]
	\arrow[from=2-1, to=2-2]
	\arrow[from=2-1, to=3-1]
	\arrow[from=2-2, to=2-3]
	\arrow[from=2-2, to=3-2]
	\arrow["{\bD (k[-1])}",from=2-3, to=3-3]
	\arrow[from=3-1, to=3-2]
	\arrow[from=3-2, to=3-3]
\end{tikzcd}\]
    upon taking cofibres vertically to obtain the bottom fibre sequence. By the bottom fibre sequence and the calculation above, we obtain a string of equivalences
\[
    \bD(X)^{(-p)}/\bD(X)^{(-p-1)}\simeq \bD(X^{(p-1)})/\bD(X^{(p)})[-1]\simeq \sideset{}{_{I_p}}\bigoplus\S[-p].
\]
    Thus, the cofibre of each map $\bD(X)^{(-p-1)}\ra \bD(X)^{(-p)}$ is equivalent to a sum of shifted spectra, and the dual CW filtration is indeed a well-defined CW filtration. Given the integral bounds $m$ and $n$ of the CW filtration of $X$, it is finite since it is clearly trivial for $p\geq n+1$, as well as for $p\leq m-1$. In the case of $p=m$, we obtain the whole of $\bD(X)$.
\end{proof}

\begin{Rmk}\label{rmk:equiv-dual-cw}
    Note that the filtration term $\bD(X)^{(-p)}$ may be equivalently expressed as follows. By \cref{lem:cofib}, we obtain an equivalence
\[
    \bD(X)^{(-p)}\coloneqq \fib(\bD(X)\ra \bD(X^{(p-1)}))\simeq \bD(\cofib(X^{(p-1)}\ra X))=\bD(X/X^{(p-1)}).
\]
    A similar argument yields an equivalence $\bD(X)^{(-p)}/\bD(X)^{(-p-1)}\simeq \bD(X^{(p-1)})/\bD(X^{(p)})$. There is slight abuse of notation---we are writing cofibres as quotients. The above identifications will be crucial in the proof of \cref{thm:A}, when defining a map between exact couples. Thus we freely interchange notationally, referring to the natural equivalence from \cref{lem:cofib}.
\end{Rmk}

\subsection*{A convenient setup for the Atiyah--Hirzebruch spectral sequence}\label{subsec:AHSS-setup}

For spectra $\cE$ and $X$, a major tool for computing the $\cE$-(co)homology groups of $X$ from the data of ordinary (co)homology is the Atiyah--Hirzebruch spectral sequence. The cohomological version of the spectral sequence is given by $E_2^{p,q}=\H^p(X;\cE^{q})\Longrightarrow \cE^{p+q}(X)$, and the homological version is $E^2_{p,q}=\H_p(X;\cE_{q})\Longrightarrow \cE_{p+q}(X)$.
There are two equivalent ways to set up the Atiyah--Hirzebruch spectral sequence. The first uses a CW filtration, the second is obtained by a Postnikov filtration of the spectrum $\cE$ as introduced in \cite{MaunderAHSS}. We opt for the former since the formulation will come in handy in \cref{sec:PD-Orientation}. A finite spectrum $X$ can be endowed with a CW filtration $0= X^{(-\infty)}\ra\dots\ra 0 \ra X^{(m)}\ra\dots \ra X^{(n)}\simeq X$ of finite spectra, and we have already seen in \cref{lem:dual-cw} how to relate it to a CW filtration on the Spanier--Whitehead dual $\bD(X)$.

We give a brief recollection of the definition of the exact couples defining the two spectral sequences. Let us fix a CW filtration $X^{(*)}$ on the finite spectrum $X$. Thus, we have extended fibration sequences (distinguished triangles) of the form
\[
    X^{(p-1)}\xlongrightarrow{i} X^{(p)}\xlongrightarrow{j}X^{(p)}/X^{(p-1)}\xlongrightarrow{k}X^{(p-1)}[1]
\]
which define the associated (unravelled) exact couple as follows. Fixing the starting degree $p+q-1$, we obtain the exact sequence
\[
    \underbrace{\cE^{p+q-1}(X^{(p)})}_{\sA_c^{p+q-1}}\xlongrightarrow{i^*}\underbrace{\cE^{p+q-1}(X^{(p-1)})}_{\sA_c^{p+q-1}}\xlongrightarrow{k[-1]^*}\underbrace{\cE^{p+q}(X^{(p)}/X^{(p-1)})}_{\sD_c^{p+q}}\xlongrightarrow{j[-1]^*}\underbrace{\cE^{p+q}(X^{(p)})}_{\sA_c^{p+q}}
\]
suppressing the suspension isomorphism $\cE^{p+q-1}((X^{(p)}/X^{(p-1)})[-1])\simeq \cE^{p+q}((X^{(p)}/X^{(p-1)}))$, and similarly $\cE^{p+q-1}(X^{(p)}[-1])\simeq \cE^{p+q}(X^{(p)})$. As a convention, we suppress suspension isomorphisms and refer to naturality of those when they appear hidden in diagrams we need to show commute. Note that the map $k[-1]^*\colon \cE^{p+q-1}(X^{(p)})\ra \cE^{p+q}(X^{(p)}/X^{(p-1)})$ is precisely the connecting homomorphism $\delta$ in the cohomological long exact sequence of the pair $(X^{(p)},X^{(p-1)})$.

We now turn to the homological exact couple. We set up the degrees to match the duality isomorphism. By \cref{lem:dual-cw}, the CW filtration on $X$ induces a dual CW filtration on its Spanier--Whitehead dual $\bD(X)$, with terms $\bD(X)^{(-p)}\coloneqq \fib(\bD(X)\ra\bD(X^{(p)}))$. We may use this filtration to define the exact couples in the homological Atiyah--Hirzebruch spectral sequence. We have extended fibration sequences (distinguished triangles) of the form
\[
    \bD(X)^{(-p-1)}\xlongrightarrow{\iota}\bD(X)^{(-p)}\xlongrightarrow{\lambda}\bD(X)^{(-p)}/\bD(X)^{(-p-1)}\xlongrightarrow{\kappa}\bD(X)^{(-p-1)}[1]
\]
which define the associated (unravelled) exact couple as follows. Fixing the starting degree $-p-q$, we obtain the exact sequence
\[
\begin{tikzcd}[column sep=0.5cm, row sep=0.2cm]
    \underbrace{\cE_{-p-q}(\bD(X)^{(-p-1)})}_{\sA^h_{-p-q}}\rar{\iota_*}\ar[draw=none]{d}[name=Y, anchor=center]{}
             & \underbrace{\cE_{-p-q}(\bD(X)^{(-p)})}_{\sA^h_{-p-q}} \rar{\lambda_*} & \underbrace{\cE_{-p-q}(\bD(X)^{(-p)}/\bD(X)^{(-p-1)})}_{\sD^h_{-p-q}} \ar["\kappa_*"', rounded corners,
                to path={ -- ([xshift=2ex]\tikztostart.east)
                	|- (Y.center) [near end]  \tikztonodes
                	-| ([xshift=-2ex]\tikztotarget.west)
                	-- (\tikztotarget)}]{dll}\\
    \underbrace{\cE_{-p-q-1}(\bD(X)^{(-p-1)})}_{\sA^h_{-p-q-1}}
               & & 
\end{tikzcd}
\]
suppressing the inverse suspension isomorphism $\cE_{-p-q}(\bD(X)^{(-p-1)}[1])\simeq \cE_{-p-q-1}(\bD(X)^{(-p-1)})$ as per our convention. The map $\smash{\kappa_*\colon \cE_{-p-q}(\bD(X)^{(-p)}/\bD(X)^{(-p-1)})\ra \cE_{-p-q-1}(\bD(X)^{(-p-1)})}$ is precisely the connecting homomorphism $\smash{\partial}$ in the homological long exact sequence of the pair $\smash{(\bD(X)^{(-p)},\bD(X)^{(-p-1)})}$.

\section{Spanier--Whitehead duality and the Atiyah--Hirzebruch spectral sequence}\label{sec:SW-AHSS}

Before we begin with the proof of \cref{thm:A}, let us give a brief outline of the strategy. The Spanier--Whitehead duality isomorphism induces an isomorphism between the groups appearing in the definition of the exact couples. This is clearly not sufficient, and we show that the Spanier--Whitehead duality isomorphism induces a morphism between the exact couples defining the spectral sequences, inducing an isomorphism between the pages. This includes the data of the differentials.

\begin{proof}[Proof of \cref{thm:A}]
    Let $X$ be a finite spectrum, together with a chosen CW filtration $0= X^{(-\infty)}\ra\dots\ra 0 \ra X^{(m)}\ra\dots \ra X^{(n)}\simeq X$ of finite spectra. By \cref{lem:dual-cw}, we obtain an induced dual CW filtration on $\bD(X)$, given by $\bD(X)^{(-p)}\coloneqq \fib(\bD(X)\ra \bD(X^{(p-1)}))$. The identification from \cref{rmk:equiv-dual-cw}, together with a simple application of \cref{lem:cofib} allows for a natural translation
\[\begin{tikzcd}[cramped,column sep=scriptsize]
	{\bD(X)^{(-p-1)}} & {\bD(X)^{(-p)}} & {\bD(X)^{(-p)}/\bD(X)^{(-p-1)}} & {\bD(X)^{(-p-1)}[1]} \\
	{\bD(X^{(p)})/\bD(X)[-1]} & {\bD(X^{(p-1)})/\bD(X)[-1]} & {\bD(X^{(p-1)})/\bD(X^{(p)})[-1]} & {\bD(X^{(p)})/\bD(X)}
	\arrow["\iota", from=1-1, to=1-2]
	\arrow["\simeq", from=1-1, to=2-1]
	\arrow["\lambda", from=1-2, to=1-3]
	\arrow["\simeq", from=1-2, to=2-2]
	\arrow["\kappa", from=1-3, to=1-4]
	\arrow["\simeq", from=1-3, to=2-3]
	\arrow["\simeq", from=1-4, to=2-4]
	\arrow["{\overline{\iota}[-1]}", from=2-1, to=2-2]
	\arrow["{\overline{\lambda}[-1]}", from=2-2, to=2-3]
	\arrow["{\overline{\kappa}[-1]}", from=2-3, to=2-4]
\end{tikzcd}\]
    of extended fibre sequences (distinguished triangles). The lower fibre sequence is induced by taking cofibrations---once more we abuse notation and write quotients instead of cofibres. The data is given by the \emph{translation diagram}
\[\begin{tikzcd}
	{\bD(X)^{(-p-1)}} & {\bD(X)} & {\bD(X^{(p)})} & {\bD(X^{(p)})/\bD(X)} \\
	{\bD(X)^{(-p)}} & {\bD(X)} & {\bD(X^{(p-1)})} & {\bD(X^{(p-1)})/\bD(X)} \\
	{\bD(X)^{(-p)}/\bD(X)^{(-p-1)}} & 0 & {\bD(X^{(p-1)})/\bD(X^{(p)})} & {\bD(X^{(p-1)})/\bD(X^{(p)})} \\
	{\bD(X)^{(-p-1)}[1]} & {\bD(X)[1]} & {\bD(X^{(p)})[1]} & {\bD(X^{(p)})/\bD(X)[1]}
	\arrow["\fib", from=1-1, to=1-2]
	\arrow["\iota", from=1-1, to=2-1]
	\arrow[from=1-2, to=1-3]
	\arrow["\id", from=1-2, to=2-2]
	\arrow["\cofib", from=1-3, to=1-4]
	\arrow["{\bD i}", from=1-3, to=2-3]
	\arrow["{\circled{1}}"{description}, draw=none, from=1-3, to=2-4]
	\arrow["{\overline{\iota}}", from=1-4, to=2-4]
	\arrow["\fib", from=2-1, to=2-2]
	\arrow["\lambda", from=2-1, to=3-1]
	\arrow[from=2-2, to=2-3]
	\arrow[from=2-2, to=3-2]
	\arrow["\cofib", from=2-3, to=2-4]
	\arrow["{\bD(k[-1])}", from=2-3, to=3-3]
	\arrow["{\circled{2}}"{description}, draw=none, from=2-3, to=3-4]
	\arrow["{\overline{\lambda}}", from=2-4, to=3-4]
	\arrow[from=3-1, to=3-2]
	\arrow["\kappa", from=3-1, to=4-1]
	\arrow[from=3-2, to=3-3]
	\arrow[from=3-2, to=4-2]
	\arrow[equals, from=3-3, to=3-4]
	\arrow["{\bD(j[-1])}", from=3-3, to=4-3]
	\arrow["{\circled{3}}"{description}, draw=none, from=3-3, to=4-4]
	\arrow["{\overline{\kappa}}", from=3-4, to=4-4]
	\arrow["{\fib[1]}", from=4-1, to=4-2]
	\arrow[from=4-2, to=4-3]
	\arrow["{\cofib[1]}", from=4-3, to=4-4]
\end{tikzcd}\]
of (co)fibre sequences. Note that the map $\bD (k[-1])$ lands in $\bD(X^{(p)}/X^{(p-1)})[1]$ which we identify with $\bD(X^{(p-1)})/\bD(X^{(p)})$, since we may write $\bD(k[-1])\simeq\bD(k)[1]$. This correspondence allows us to equivalently define the homological exact couple using the right-most (shifted) cofibre sequence. This allows for a more natural translation between the cohomological exact couple $(\sA_c^{*-1},\sD_c^{*})$ and homological exact couples $(\sA^h_{*},\sD^h_*)$. We number several of the squares in the translation diagram for future referencing. For a couple map $(\Psi, \Phi)\colon(\sA_c^{*-1},\sD_c^*)\ra (\sA^h_{-*}, \sD^h_{-*})$ to induce an isomorphism between the two associated spectral sequences, it suffices that $\Phi$ is an isomorphism. We define
\[
    \Phi\colon \cE^{p+q}(X^{(p)}/X^{(p-1)})\xrightarrow{\SW}\cE_{-p-q}(\bD(X^{(p)}/X^{(p-1)}))\xlongrightarrow{\cref{lem:cofib}}\cE_{-p-q+1}(\bD(X^{(p-1)})/\bD(X^{(p)}))
\]
where we, by convention, suppress the suspension isomorphism. Note that we have a natural identification $\cE_{-p-q+1}(\bD(X^{(p-1)})/\bD(X^{(p)}))\simeq\cE_{-p-q}(\bD(X)^{(-p)}/\bD(X)^{(-p-1)})$, and thus the composition defines, by slight abuse of notation, a map $\Phi\colon \sD_c^*\ra \sD^h_{-*}$. Furthermore, by the calculation from \cref{lem:dual-cw}, we compute both groups in question via (de)suspending accordingly, and we may interpret $\Phi$ as the Spanier--Whitehead duality isomorphism $\SW\colon \bigoplus_{I_p}\cE^q\ra\bigoplus_{I_p}\cE_{-q}$.
Furthermore, we define
\[
    \Psi\colon\cE^{p+q-1}(X^{(p-1)})\xrightarrow{\SW} \cE_{-p-q+1}(\bD(X^{(p-1)}))\xrightarrow{\cofib_*}\cE_{-p-q+1}(\bD(X^{(p-1)})/\bD(X)).
\]
Once again, there is a natural identification $\cE_{-p-q+1}(\bD(X^{(p-1)})/\bD(X))\simeq \cE_{-p-q}(\bD(X)^{(-p)})$, in which we suppress an inverse suspension isomorphism. This is mirroring the suppressed suspension isomorphism in the definition of $\Phi$. Thus, the composition induces, by abuse of notation, a map $\Psi\colon \sA_c^{*-1}\ra \sA^h_{-*}$. It is left to show that the pair $(\Psi, \Phi)\colon(\sA_c^{*-1},\sD_c^*)\ra (\sA^h_{-*}, \sD^h_{-*})$ defines a couple map. To do so, we need to show commutativity of three diagrams. The first diagram
\[\begin{tikzcd}
	{\cE^{p+q-1}(X^{(p)})} & {\cE^{p+q-1}(X^{(p-1)})} \\
	{\cE_{-p-q+1}(\bD(X^{(p)}))} & {\cE_{-p-q+1}(\bD(X^{(p-1)}))} \\
	{\cE_{-p-q+1}(\bD(X^{(p)})/\bD(X))} & {\cE_{-p-q+1}(\bD(X^{(p-1)})/\bD(X))} 
	\arrow["{i^*}", from=1-1, to=1-2]
	\arrow["\SW", from=1-1, to=2-1]
	\arrow["\SW", from=1-2, to=2-2]
	\arrow["{\bD i_*}", from=2-1, to=2-2]
	\arrow["{\cofib_*}", from=2-1, to=3-1]
	\arrow["{\cofib_*}", from=2-2, to=3-2]
	\arrow["{\overline{\iota}_*}",from=3-1, to=3-2]
\end{tikzcd}\]
commutes by naturality of the Spanier--Whitehead duality isomorphism and commutativity of the upper-right square $\scaleobj{0.9}{\circled{1}}$ in the translation diagram above. The second square
\[\begin{tikzcd}
	{\cE^{p+q-1}(X^{(p-1)})} & {\cE^{p+q}(X^{(p)}/X^{(p-1)})} \\
	{\cE_{-p-q+1}(\bD(X^{(p-1)}))} & {\cE_{-p-q}(\bD(X^{(p)}/X^{(p-1)}))} \\
	{\cE_{-p-q+1}(\bD(X^{(p-1)})/\bD(X))} & {\cE_{-p-q+1}(\bD(X^{(p-1)})/\bD(X^{(p)}))} 
	\arrow["{k[-1]^*}", from=1-1, to=1-2]
	\arrow["\SW", from=1-1, to=2-1]
	\arrow["\SW", from=1-2, to=2-2]
	\arrow["{\bD (k[-1])_*}", from=2-1, to=2-2]
	\arrow["{\cofib_*}", from=2-1, to=3-1]
	\arrow["{\simeq}",from=2-2, to=3-2]
	\arrow["{\overline{\lambda}_*}",from=3-1, to=3-2]
\end{tikzcd}\]
commutes as well by naturality of the Spanier--Whitehead duality isomorphism and commutativity of the middle-right square $\scaleobj{0.9}{\circled{2}}$ in the translation diagram above. The last diagram 
\[\begin{tikzcd}
	{\cE^{p+q}(X^{(p)}/X^{(p-1)})} & {\cE^{p+q}(X^{(p)})} \\
	{\cE_{-p-q}(\bD(X^{(p)}/X^{(p-1)}))} & {\cE_{-p-q}(\bD(X^{(p)}))} \\
	{\cE_{-p-q+1}(\bD(X^{(p-1)})/\bD(X^{(p)}))} & {\cE_{-p-q}(\bD(X^{(p-1)})/\bD(X))}
	\arrow["{j[-1]^*}", from=1-1, to=1-2]
	\arrow["\SW", from=1-1, to=2-1]
	\arrow["\SW", from=1-2, to=2-2]
	\arrow["{\bD(j[-1])_*}", from=2-1, to=2-2]
	\arrow["\simeq", from=2-1, to=3-1]
	\arrow["{\cofib[1]_*}", from=2-2, to=3-2]
	\arrow["{\overline{\kappa}_*}", from=3-1, to=3-2]
\end{tikzcd}\]
commutes by naturality of the Spanier--Whitehead duality isomorphism and commutativity of the diagram lower-right square $\scaleobj{0.9}{\circled{3}}$ in the translation diagram. This concludes the argument that the pair $(\Psi, \Phi)\colon(\sA_c^{*-1},\sD_c^*)\ra (\sA^h_{-*}, \sD^h_{-*})$ defines a couple map. The couple map induces an isomorphism of spectral sequences since $\Phi\colon \sD_c^*\ra \sD^h_{-*}$ is an isomorphism. For completion, the diagram
\[\begin{tikzcd}
	{C_{\text{cell}}^{p-1}(X;\cE^q)} && {C_{\text{cell}}^{p}(X;\cE^q)} \\
	{\bigoplus_{I_{p-1}}\wt{\cE}^{p-1+q}(\S[{p-1]})} && {\prod_{I_{p}}\wt{\cE}^{p+q}(\S[p])} \\
	{\cE^{p+q-1}(X^{(p-1)}/X^{(p-2)})} & {\cE^{p+q-1}(X^{(p-1)})} & {\cE^{p+q}(X^{(p)}/X^{(p-1)})} \\
	{\cE_{-p-q+1}(\bD^{(-p-1)}/\bD^{(-p-2)})} & {\cE_{-p-q+1}(\bD^{(-p-1)})} & {\cE_{-p-q}(\bD^{(-p)}/\bD^{(-p-1)})} \\
	{\prod_{I_{p-1}}\wt{\cE}_{-p-q+1}(\S[{-p+1]})} && {\prod_{I_{p}}\wt{\cE}_{-p-q}(\S[{-p]})} \\
	{C^{\text{cell}}_{-p+1}(\bD(X);\cE_{-q})} && {C^{\text{cell}}_{-p}(\bD(X);\cE_{-q})}
	\arrow["{d_1}", dashed, from=1-1, to=1-3]
	\arrow["\simeq", from=1-1, to=2-1]
	\arrow["\simeq", from=1-3, to=2-3]
	\arrow["\simeq", from=2-1, to=3-1]
	\arrow["\simeq", from=2-3, to=3-3]
	\arrow["{j[-1]^*}", from=3-1, to=3-2]
	\arrow["\Phi", from=3-1, to=4-1]
	\arrow["{k[-1]^*}", from=3-2, to=3-3]
	\arrow["\Psi", from=3-2, to=4-2]
	\arrow["\Phi", from=3-3, to=4-3]
	\arrow["{\kappa_*}", from=4-1, to=4-2]
	\arrow["{\lambda_*}", from=4-2, to=4-3]
	\arrow["\simeq"', from=5-1, to=4-1]
	\arrow["\simeq"', from=5-3, to=4-3]
	\arrow["\simeq"', from=6-1, to=5-1]
	\arrow["{d^1}", dashed, from=6-1, to=6-3]
	\arrow["\simeq"', from=6-3, to=5-3]
\end{tikzcd}\]
depicts the translation between the first page of both spectral sequences.
\end{proof}

\begin{Rmk}\label{rmk:r-w}
    As pointed out by Oscar Randal-Williams, the argument can be drastically simplified by using the definition of the Atiyah--Hirzebruch spectral sequences via a Postnikov filtration of the spectrum $\cE$ \cite{MaunderAHSS}. Maunder only defines the cohomological Atiyah--Hirzebruch spectral sequence in such a way, though one can define the homological Atiyah--Hirzebruch spectral sequence in a similar manner---essentially because the smash product $-\otimes-\colon \Sp\times \Sp\ra \Sp$ preserves cofibrations. With this setup, one may directly identify the exact couples via the adjunction inducing the Spanier--Whitehead duality isomorphism by replacing $\cE$-(co)homology by the (shifted) Postnikov truncations, namely $\tau_{\leq n}(\cE)$-(co)homology.
\end{Rmk}

\section{Orientations on Poincar\'e duality complexes}\label{sec:PD-Orientation}

\subsection*{Orientations on manifolds}

To motivate the definition of an $\cR$-orientation on a Poincar\'e duality complex $X$ for a ring spectrum $\cR$, we begin with a brief discussion on orientations of vector bundles over CW complexes, and then focus on manifolds. We refer the reader to \cite[{Chapter~V}]{RudyakThomSpectra} for details. This is followed by a survey on generalisations to spherical fibrations, and Poincar\'e duality complexes.

\begin{Def}\label{def:E-orientation}
Let $\cR$ be a ring spectrum, $X$ a CW complex together with a vector bundle $\xi\colon X\ra \BO(r)$ of rank $r$. The bundle $\xi$ is \emph{$\cR$-orientable} if a factorisation of the unit map $\mu\colon \S\ra \cR$ of the form
\[\begin{tikzcd}
	\S & \cR \\
	{\Sigma^\infty\S^{r}[-r]} & {\Sigma^\infty\Th(\xi)[-r]}
	\arrow["\mu", from=1-1, to=1-2]
	\arrow["\simeq", from=1-1, to=2-1]
	\arrow["{\iota_*}", from=2-1, to=2-2]
	\arrow["{u(\xi)}"', from=2-2, to=1-2]
\end{tikzcd}\]
    exists. The map $\iota_*$ is induced from the canonical fibre inclusion $\iota\colon \S^r\ra \Th(\xi)$. The map $u(\xi)\colon \Sigma^\infty\Th(\xi)[-r]\ra \cR$ defines a \emph{Thom class} $u(\xi)\in \cR^r(\Sigma^\infty\Th(\xi))$ via the chain of equivalences
\begin{align*}
    u(\xi)\in [\Sigma^\infty\Th(\xi)[-r],\cR] & \simeq [\Sigma^\infty\Th(\xi)\otimes \S[-r],\cR] \simeq [\S[-r],\F(\Sigma^\infty\Th(\xi),\cR)] \\
    & \simeq \pi_{-r}(\F(\Sigma^\infty\Th(\xi),\cR))=\cR^{r}(\Sigma^\infty\Th(\xi)).
\end{align*}
A choice of such a Thom class is an \emph{$\cR$-orientation} of the bundle $\xi$.
\end{Def}

It is a classical result that a Thom class $u(\xi)\in \cR^r(\Sigma^\infty\Th(\xi))\simeq \cR^r(\bD(\xi),\S(\xi))$ gives rise to a Thom isomorphism as follows. Let $\pi\colon \bD(\xi)\ra X$ be the canonical projection from the disc bundle of $\xi$. Then the composition
\[
    \cR^n(X)\xlongrightarrow{\pi^*}\cR^n(\Sigma^\infty\bD(\xi)_+)\xlongrightarrow{-\cup u(\xi)}\cR^{n+r}(\Sigma^\infty\Th(\xi))
\]
is an isomorphism.
For a closed, $d$-dimensional manifold $M$, we say that $M$ is \emph{$\cR$-orientable} if such a factorisation as above exists for (a representative of) the stable normal bundle $\nu_M\colon M\ra \BO$. In this case, the data suffices to define an $\cR$-fundamental class $[M]_\cR\in \cR_d(M)$ as follows.

\smallskip
\noindent\textit{Orientation and fundamental class.}
Let $i\colon M\hra \R^{\ell}$ be an embedding of a closed, $d$-manifold $M$ into high-dimensional Euclidean space, such that its associated normal bundle $\nu_i\colon M\ra \BO(\ell-d)$ represents the stable normal bundle $\nu_M\colon M\ra \BO$ as a virtual bundle.
An $\cR$-orientation on $M$ is equivalent to a Thom class $u(\nu_i)$ in $\cR^{\ell-d}(\Sigma^\infty\Th(\nu_i))$ which has the defining property that the restriction along the fibre inclusion $\iota\colon \S^{\ell-d}\ra \Th(\nu_i)$ sends the class $u(\nu_i)$ to the unit of the ring $\wt{\cR}^{\ell-d}(\S^{\ell-d})\simeq \pi_0(\cR)$.
Since the Spanier--Whitehead dual of $\Sigma^\infty\Th(\nu_i)[-\ell]$ is precisely $\Sigma^\infty M_+$ (\cite[{Proposition~3.2}]{AtiyahThomComplexes}), which is classically known as \emph{Atiyah duality}, we obtain an associated homology class as the image of $u(\nu_i)$ under the equivalence $\cR^{\ell-d}(\Sigma^\infty \Th(\nu_i))\simeq \cR_d(M)$ induced by adjunction properties of Spanier--Whitehead duality---namely, the equivalences
\begin{align*}
    \cR^{\ell-d}(\Sigma^\infty\Th(\nu_i)) & \simeq [\Sigma^\infty\Th(\nu_i)[-(\ell-d)],\cR] \simeq [\Sigma^\infty\Th(\nu_i)[-\ell]\otimes \S[d],\cR] \\
    & \simeq [\S[d],\cR\otimes \Sigma^\infty M_+]\simeq\cR_d(M).
\end{align*}
This class does not depend on the embedding $i\colon M\hra \R^{\ell}$ since any two embeddings are isotopic for large enough $\ell$, and thus the Thom spaces of the normal bundles are equivalent. Alternatively, one may pass to the virtual bundle $-TM$ via a split of the exact sequence
\[
    TM\lra T\R^{\ell}|_M\lra \nu_i
\]
of vector bundles. We call the homological class the \emph{$\cR$-fundamental class} $[M]_\cR\in \cR_d(M)$.

\smallskip
\noindent\textit{Poincar\'e duality.}
The multiplicative structure on the ring spectrum induces a cap product pairing, and the Poincar\'e duality isomorphism
\[
    \PD\colon\cR^n(M)\xlongrightarrow{-\cap [M]_\cR}\cR_{d-n}(M).
\]
is obtained by sending a cohomology class $x\in \cR^n(M)$ to the homology class $x\cap [M]_\cR\in \cR_{d-n}(M)$ as shown in \cite[{Chapter~V,~Theorem~2.9}]{RudyakThomSpectra}. In general, the argument can be adjusted to allow manifolds with non-empty boundary. In this case, the fundamental class $[M,\partial M]_\cR$ lies in $\cR_d(M,\partial M)$, and the Poincar\'e-Lefschetz duality isomorphism $\PD\colon \cR^n(M,\partial M)\ra \cR_{d-n}(M)$ is obtained by forming the cap-product with the relative fundamental class $[M,\partial M]_\cR$.

The way we obtain the $\cR$-fundamental class $[M]_\cR$ from the Thom class $u(\nu_i)$ together with an application of the Spanier--Whitehead duality isomorphism, leads to the factorisation
\[
    \PD\colon \cR^n(M)\xlongrightarrow{\text{Thom}}\cR^{n+l-d}(\Sigma^\infty\Th(\nu_i))\simeq\cR^{n-d}(\Sigma^\infty\Th(\nu_i)[-\ell])\xlongrightarrow{\SW}\cR_{d-n}(M)
\]
of the Poincar\'e duality isomorphism.

\begin{Ex}
    We provide a short collection of important notions and results.
    \begin{enumerate}
        \item[(1)] The classical notion of a \emph{complex oriented cohomology theory} $\cE$ ensures that every complex vector bundle admits an $\cE$-orientation, the universal complex oriented cohomology theory being famously given by $\MU$.
        \item[(2)] A \emph{framed manifold} $M$ is $\cR$-orientable for any ring spectrum $\cR$. A factorisation of the unit map $\mu\colon\S\ra \cR$ can be constructed as follows. A trivialisation of the normal bundle associated to the embedding $i\colon M\hra \R^{d+\ell}$ induces an equivalence $\Sigma^\infty \Th(\nu_i)[-\ell]\simeq \Sigma^\infty M_+$. The fibre inclusion map $\iota_*\colon \Sigma^\infty\S^\ell[-\ell]\ra \Sigma^\infty \Th(\nu_i)[-\ell]\simeq \Sigma^\infty M_+$ admits a left split via the collapse map $\Sigma^\infty M_+\ra \S$. As an orientation, one can choose the composition of the collapse map followed by the unit map. The same argument can be done for a manifold $M$ with \emph{trivial Spivak normal fibration} since the associated Thom space admits the description as above.
        \item[(3)] For cobordism theories, there is a tower of orientations given by the Whitehead filtration of the ``Thomification'' of $\BO$
\[\begin{tikzcd}
	\S & \MFr & \dots & \MString & \MSpin & \MSO & \MO \\
	&&& \tmf & \ko & {\H\Z} & {\H\Z/2}
	\arrow[equals, from=1-1, to=1-2]
	\arrow[from=1-2, to=1-3]
	\arrow[from=1-3, to=1-4]
	\arrow[from=1-4, to=1-5]
	\arrow["\sigma", from=1-4, to=2-4]
	\arrow[from=1-5, to=1-6]
	\arrow["{\widehat{A}}", from=1-5, to=2-5]
	\arrow[from=1-6, to=1-7]
	\arrow[from=1-6, to=2-6]
	\arrow[from=1-7, to=2-7]
\end{tikzcd}\]
        with classical orientations on $\MO$ and $\MSO$. The equivalence $\S\simeq \MFr$ obtained from the Pontryagin--Thom construction yields the previous point via the unit map $\mu\colon \S\ra \cR$ for any ring spectrum $\cR$. The $7$-connected map $\smash{\widehat{A}}\colon\MSpin\ra \ko$ represents the Atiyah--Bott--Shapiro $\smash{\widehat{A}}$-orientation as developed in \cite{ABS}, and \cite{JoachimThesis}. Furthermore, this map of spectra refines to a map of $\E_\infty$-ring spectra as shown in \cite{JoachimKTheory} (even equivariantly). The orientation $\sigma\colon \MString\ra \tmf$, which induces the Witten genus on homotopy groups, is due to work of Ando--Hopkins--Rezk--Strickland (see for example \cite{AHR}, \cite{AHS}). Thus, spin manifolds are $\ko$-orientable, and string manifolds admit a $\tmf$-orientation.
        \item[(4)] By \cite[{Chapter~16}]{RanickiLTheory}, every manifold that is classically orientable (orientable over $\H\Z$) is orientable over the connective symmetric $\L$-theory spectrum $\LL^s(\Z)\langle 0\rangle$ of the integers. This is expressed via an appropriate map $\sigma_{\LL}\colon\MSTop\ra \LL^s(\Z)\langle 0\rangle$.
    \end{enumerate}
\end{Ex}

\subsection*{Passing to spherical fibrations and Poincar\'e duality complexes}

A Poincar\'e duality complex $X$ does not have a stable normal bundle in general. Its replacement is the Spivak normal fibration $\SF(X)\colon X\ra \BG$, which is a stable spherical fibration. Since spherical fibrations have well-defined Thom spaces, we may extend \cref{def:E-orientation} to spherical fibrations as follows.

\begin{Def}\label{def:spherical-orientation}
    Let $\cR$ be a ring spectrum, $X$ a CW complex together with a spherical fibration $\zeta\colon X\ra \BG(r)$ of rank $r$. The spherical fibration $\zeta$ is \emph{$\cR$-orientable} if a factorisation of the unit map $\mu\colon \S\ra \cR$ of the form
\[\begin{tikzcd}
	\S & \cR \\
	{\Sigma^\infty\S^{r}[-r]} & {\Sigma^\infty\Th(\zeta)[-r]}
	\arrow["\mu", from=1-1, to=1-2]
	\arrow["\simeq", from=1-1, to=2-1]
	\arrow["{\iota_*}", from=2-1, to=2-2]
	\arrow["{u(\zeta)}"', from=2-2, to=1-2]
\end{tikzcd}\]
    exists. The map $\iota_*$ is induced from the canonical fibre inclusion $\iota\colon \S^r\ra \Th(\zeta)$. The map $u(\zeta)\colon \Sigma^\infty\Th(\zeta)[-r]\ra \cR$ defines a \emph{Thom class} $u(\zeta)\in \cR^r(\Sigma^\infty\Th(\zeta))$ as in \cref{def:E-orientation}. A choice of such a Thom class is an \emph{$\cR$-orientation} of the spherical fibration $\zeta$.
\end{Def}

As before, such a Thom class $u(\zeta)\in \cR^r(\Sigma^\infty\Th(\zeta))$ gives rise to a Thom isomorphism
\[
    \cR^n(X)\xlongrightarrow{-\cup u(\zeta)}\cR^{n+r}(\Sigma^\infty \Th(\zeta)),
\]
via the map $\Delta^r\colon\Sigma^\infty\Th(\zeta)\rightarrow \Sigma^\infty X_+\otimes \Sigma^\infty\Th(\zeta)$ which arises from the space-level identification $\Th(\zeta\times\varepsilon^1)\simeq \Th(\zeta)\wedge X_+$. In particular, there is a pairing
\[
    \cR^n(\Sigma^\infty X_+)\otimes\cR^r(\Sigma^\infty\Th(\zeta))\xlongrightarrow{-\cup-} \cR^{n+r}(\Sigma^\infty X_+\otimes \Sigma^\infty\Th(\zeta))\xlongrightarrow{(\Delta^r)^*}\cR^{n+r}(\Sigma^\infty\Th(\zeta)).
\]
The Thom isomorphism takes as input an arbitrary class $x\in \cR^n(\Sigma^\infty X_+)$, and the Thom class $u(\zeta)\in \cR^r(\Sigma^\infty\Th(\zeta))$. For details, see for example \cite[{Chapter~V,~Theorem-Definition~1.3}]{RudyakThomSpectra}, or \cite[{Theorem~20.5.8}]{ParametrizedHomotopyTheory}. Furthermore, we may express the Thom isomorphism as an equivalence
\[
    \TD_h\colon\cR\otimes \Sigma^\infty\Th(\zeta)\lra\cR\otimes \Sigma^\infty X_+[r]
\]
of spectra \cite{ANoteOnThom}; evidently, there is also an equivalence of mapping spectra
\[
\TD_c\colon\F(\Sigma^\infty X_+,\cR)\lra \F(\Sigma^\infty\Th(\zeta)[-r],\cR).
\]
These equivalences induce the classical Thom isomorphism upon taking homotopy groups. We call the map $\TD_h$ the \emph{homological Thom equivalence}, and the map $\TD_c$ the \emph{cohomological Thom equivalence}---note that one induces the other, and thus we focus on the cohomological Thom equivalence.

The analogue of the stable normal bundle of a manifold for a Poincar\'e duality complex $X$ is the Spivak normal fibration $\SF(X)\colon X\ra \BG$. Hence, we say a Poincar\'e duality complex $X$ is \emph{orientable over a ring spectrum $\cR$} if (a representative of) its Spivak normal fibration is orientable in the sense of \cref{def:spherical-orientation}. 

In the following, we unravel the definition to obtain a Thom isomorphism and equivalence for the \emph{stable} Spivak normal fibration. We carefully go through the construction, as there are different conventions in the literature.

\smallskip
\noindent\textit{Stable Thom isomorphisms.}
We briefly discussed the construction of the Spivak normal fibration in \cref{sec:preliminaries}. Let us consider an embedding $i\colon X\hra \R^\ell$ of a Poincar\'e duality complex $X$ into high-dimensional Euclidean space. There is an associated spherical fibration $\SF_\ell(X)\colon X\ra \BG(\ell-d)$ which is a representative of the Spivak normal fibration $\SF\colon X\ra \BG$ in the group $[X,\BG]\simeq \Pic(\S)^0(X)$. Note that we suppress the embedding in the notation of the spherical fibration.

\begin{Def}\label{def:thom-SF}
    Let $X$ be a $d$-dimensional Poincar\'e duality complex with Spivak normal fibration $\SF(X)\colon X\ra \BG$. The Thom spectrum of the Spivak normal fibration is defined as
\[
    \MSF(X)\coloneqq\colim_{N\geq\ell}\Sigma^\infty\Th(\SF_N(X))[-N]
\]
    where we have maps $\Th(\SF_\ell(X))\ra \Th(\SF_{\ell+1}(X))\ra\dots$ induced by inclusions into higher-dimensional Euclidean space.
\end{Def}

Since we assume $\ell$ to be large enough that $\SF_\ell(X)$ represents the Spivak normal fibration, we immediately obtain the equivalence $\SF_{\ell+1}(X)\simeq \SF_{\ell}\oplus\varepsilon^1$. On Thom spaces, this corresponds to suspension---namely, we have an equivalence $\Th(\SF_{\ell+1}(X))\simeq \Sigma\Th(\SF_{\ell}(X))$, and on suspension spectra $\Sigma^\infty \Th(\SF_{\ell+1}(X))\simeq \Sigma^\infty\Th(\SF_{\ell}(X))[1]$. Thus, the procedure for obtaining the spectrum $\MSF(X)$ is the same as \emph{spectrification} of the prespectrum given by the sequence of spaces $\{\Th(\SF_{N}(X))\}_{N\geq \ell}$ with appropriate structure maps.

\begin{Lem}\label{lem:MSF-orientation}
    Let $\cR$ be a ring spectrum, and $X$ a $d$-dimensional Poincar\'e duality complex together with its Spivak normal fibration $\SF(X)\colon X\ra \BG$. An $\cR$-orientation on (a representative of) the Spivak normal fibration induces a Thom class $u(\SF(X))\in \cR^{-d}(\MSF(X))$ and a Thom equivalence
\[
    \TD_c\colon \F(\Sigma^\infty X_+,\cR)\ra \F(\MSF(X)[d],\cR)
\]
    inducing the classical Thom isomorphism $-\cup u(\SF(X))\colon\cR^{n}(X)\ra\cR^{n-d}(\MSF(X))$ on homotopy groups.
\end{Lem}

\begin{proof}
    As usual, we begin with an embedding $i\colon X\hra \R^\ell$ of $X$ into high-dimensional Euclidean space, which yields a spherical fibration $\SF_\ell(X)\colon X\ra \BG(\ell-d)$ representing the Spivak normal fibration. By assumption, this spherical fibration is orientable over $\cR$ in the sense of \cref{def:spherical-orientation}, and we may choose a Thom class $u(\SF_\ell(X))$. This Thom class allows for a Thom isomorphism
\[
    -\cup u(\SF_\ell(X))\colon \cR^n(X)\lra\cR^{n+\ell-d}(\Sigma^\infty \Th(\SF_\ell(X))).
\]
    The Thom class $u(\SF(X))\in \cR^{-d}(\MSF(X))$ is obtained as follows. Orientations are compatible with stabilisations (\cite[Chapter~V,~Proposition~1.10]{RudyakThomSpectra}), and thus we have a sequence
\[
    \dots\longrightarrow\cR^{\ell+1-d}(\Sigma^\infty\Th(\SF_{\ell+1}(X)))\longrightarrow \cR^{\ell-d}(\Sigma^\infty\Th(\SF_{\ell}(X)))
\]
    mapping Thom classes to Thom classes. Taking the limit of the sequence, we obtain a class
\begin{align*}
    u(\SF(X)) &\in \lim( \dots\longrightarrow\cR^{\ell+1-d}(\Sigma^\infty\Th(\SF_{\ell+1}(X)))\longrightarrow \cR^{\ell-d}(\Sigma^\infty\Th(\SF_{\ell}(X))))\\
    & \simeq \cR^{-d}(\colim(\Sigma^\infty\Th(\SF_{\ell}(X))[-\ell]\longrightarrow \Sigma^\infty\Th(\SF_{\ell+1}(X))[-(\ell+1)]\longrightarrow\dots))\\
    & =\cR^{-d}(\MSF(X))
\end{align*}
    which is the induced Thom class of the Spivak normal fibration of $X$. Similarly, we have a system of Thom equivalences
\[\begin{tikzcd}
	\vdots & \vdots \\
	{\F(\Sigma^\infty X_+,\cR)} & {\F(\Sigma^\infty \Th(\SF_{\ell+1}(X))[-(\ell+1-d)],\cR)} \\
	{\F(\Sigma^\infty X_+,\cR)} & {\F(\Sigma^\infty \Th(\SF_\ell(X))[-(\ell-d)],\cR)}
	\arrow[from=1-2, to=2-2]
	\arrow[equals, from=2-1, to=1-1]
	\arrow["{\TD_c}", from=2-1, to=2-2]
	\arrow[from=2-2, to=3-2]
	\arrow[equals, from=3-1, to=2-1]
	\arrow["{\TD_c}", from=3-1, to=3-2]
\end{tikzcd}\]
which, upon taking vertical limits, induces an equivalence
\begin{align*}
    \TD_c\colon \F(\Sigma^\infty X_+,\cR) &\lra \lim_{N\geq \ell}(\F(\Sigma^\infty\Th(\SF_N(X))[-(N-d)],\cR))\\
    & \simeq \F(\colim_{N\geq \ell}(\Sigma^\infty\Th(\SF_N(X))[-(N-d)]),\cR)\\
    & \simeq \F(\MSF(X)[d],\cR).
\end{align*}
Notice that we almost repeated the same argument twice---the existence of the Thom class, and the Thom equivalence are shown in the same way. By \cite{ANoteOnThom}, this yields the classical Thom isomorphism $-\cup u(\SF(X))\colon\cR^{n}(X)\ra\cR^{n-d}(\MSF(X))$.
\end{proof}

\subsection*{Poincar\'e duality and the Atiyah--Hirzebruch spectral sequence}\label{subsec:PD-AHSS}

This subsection contains the proof of \cref{thm:B}. The strategy is to factor Poincar\'e duality as the Thom isomorphism, followed by Spanier--Whitehead duality. By \cref{thm:A}, we know that the latter induces an isomorphism between the Atiyah--Hirzebruch spectral sequences, and it is left to show that the same holds for the Thom isomorphism. Before we may identify the exact couples, we begin with a well-known fact on restricting the Thom isomorphism to skeleta.

\begin{Lem}\label{lem:skeleta-orientations}
    Let $X$ be a CW complex, and $\zeta\colon X\ra \BG(r)$ be a spherical fibration that is oriented over a ring spectrum $\cR$. The restriction $\zeta^{(p)}\coloneqq\zeta|_{X^{(p)}}\colon X^{(p)}\ra \BG(r)$ to the $p$-skeleton $X^{(p)}$ of $X$ inherits an $\cR$-orientation, and the Thom spectrum $\Sigma^\infty\Th(\zeta^{(p)})$ agrees with $\Sigma^\infty \Th(\zeta)^{(p+r)}$. Furthermore, there is a Thom equivalence
\[
    \TD_c\colon \F(\Sigma^\infty X^{(p)}_+,\cR)\lra \F(\Sigma^\infty \Th(\zeta)^{(p+r)}[-r],\cR)
\]
    of spectra which induces the restricted Thom isomorphism on homotopy groups.
\end{Lem}

\begin{proof}
    The fact that the restricted spherical fibration $\zeta^{(p)}\coloneqq\zeta|_{X^{(p)}}\colon X^{(p)}\ra \BG(r)$ inherits an $\cR$-orientation immediately follows from \cite[Chapter~V,~Proposition~1.10]{RudyakThomSpectra} applied to the inclusion $X^{(p)}\hra X$. The statement about the Thom spectra is already true on spaces---the Thom space $\Th(\zeta^{(p)})$ is equivalent to the $(p+r)$-skeleton of the Thom space $\Th(\zeta)$. This is a classical result for vector bundles, and it may be shown by considering the open $p$-cells on $X$, together with the fact that a vector bundle is locally trivial. The vector space structure on the fibre is not used in the proof, and thus we may replace the r\^{o}le of the vector bundle by the disc bundle $D(\zeta)$ of the sphere bundle $\zeta$, defined as the fibre-wise cone. Thus, one obtains a commutative diagram of the form
\[\begin{tikzcd}
	{\D^p\times \S^{r-1}} & {S(i^*\zeta)} & {S(\zeta)} \\
	{\D^p\times \D^{r}} & {D(i^*\zeta)} & {D(\zeta)} & {\Th(\zeta)} \\
	& {\D^p} & X
	\arrow["\simeq", from=1-1, to=1-2]
	\arrow[hook, from=1-1, to=2-1]
	\arrow[from=1-2, to=1-3]
	\arrow[hook, from=1-2, to=2-2]
	\arrow[hook, from=1-3, to=2-3]
	\arrow["\simeq", from=2-1, to=2-2]
	\arrow[from=2-2, to=2-3]
	\arrow[from=2-2, to=3-2]
	\arrow[from=2-3, to=2-4]
	\arrow[from=2-3, to=3-3]
	\arrow["i", from=3-2, to=3-3]
\end{tikzcd}\]
    where the outer square in the middle is associated to the pullback of the spherical fibration $\zeta$ under characteristic map $i\colon \D^p\ra X$ of a $p$-cell.
\end{proof}

Let $X$ be a $d$-dimensional Poincar\'e duality complex that is oriented over a ring spectrum $\cR$, and thus comes with a Thom class $u(\SF_\ell(X))\in \cR^{\ell-d}(\Sigma^\infty \Th(\SF_\ell(X)))$. Consider the (unstable) Thom isomorphism $-\cup u(\SF_\ell(X))\colon\cR^n(X)\ra\cR^{n+\ell-d}(\Sigma^\infty \Th(\SF_\ell(X)))$ as in the proof of \cref{lem:MSF-orientation}. We begin with a morphism between the Atiyah--Hirzebruch spectral sequences induced by the unstable Thom isomorphism, and then pass to a limit argument to obtain the result for the stable Thom isomorphism. Recall the setup of the Atiyah--Hirzebruch spectral sequence from \cref{subsec:AHSS-setup}. The (unravelled) exact couple defining the cohomological Atiyah--Hirzebruch spectral sequence computing $\cR^*(X)$ yields an exact sequence
\[
    \underbrace{\cR^{p+q-1}(X^{(p)})}_{\sA_X^{p+q-1}}\xlongrightarrow{i^*}\underbrace{\cR^{p+q-1}(X^{(p-1)})}_{\sA_X^{p+q-1}}\xlongrightarrow{k[-1]^*}\underbrace{\cR^{p+q}(X^{(p)}/X^{(p-1)})}_{\sD_X^{p+q}}\xlongrightarrow{j[-1]^*}\underbrace{\cR^{p+q}(X^{(p)})}_{\sA_X^{p+q}}
\]
upon fixing the starting degree $p+q-1$.

Writing $\wt{\cR}^*(\Th(\SF_\ell(X)))= \cR^*(\Sigma^\infty \Th(\SF_\ell(X)))$ for brevity, the exact couple defining the Atiyah--Hirzebruch spectral sequence computing $\cR^*(\Sigma^\infty \Th(\SF_\ell)(X))$ yields an exact sequence
\[
\begin{tikzcd}[column sep=tiny, row sep=0.2cm]
    \underbrace{\wt{\cR}^{p+q-1+\ell-d}(\Th(\SF_\ell(X))^{(p+\ell-d)})}_{\sA_{\SF_\ell(X)}^{p+q-1+\ell-d}}\rar{\iota^*}\ar[draw=none]{d}[name=Y, anchor=center]{}
             & \underbrace{\wt{\cR}^{p+q+\ell-d}(\Th(\SF_\ell(X))^{(p-1+\ell-d)})}_{\sA_{\SF_\ell(X)}^{p+q+\ell-d}} 
             \ar["{\kappa[-1]^*}"', rounded corners,
                to path={ -- ([xshift=2ex]\tikztostart.east)
                	|- (Y.center) [near end]  \tikztonodes
                	-| ([xshift=-2ex]\tikztotarget.west)
                	-- (\tikztotarget)}]{dl}\\  
    \underbrace{\wt{\cR}^{p+q+\ell-d}(\Th(\SF_\ell(X))^{(p+\ell-d)}/\Th(\SF_\ell(X))^{(p-1+\ell-d)})}_{\sD_{\SF_\ell(X)}^{p+q+\ell-d}} \rar{\lambda[-1]^*} &
    \underbrace{\wt{\cR}^{p+q+\ell-d}(\Th(\SF_\ell(X))^{(p+\ell-d)})}_{\sA_{\SF_\ell(X)}^{p+q+\ell-d}}
\end{tikzcd}
\]
upon fixing the starting degree $p+q-1+\ell-d$.

\begin{proof}[{Proof of \cref{thm:B}}]
The setup from \cref{sec:PD-Orientation} allows us to express Poincar\'e duality as follows. Let $X$ be a Poincar\'e duality complex of dimension $d$, oriented over a ring spectrum $\cR$. By \cite[{Remark~A.8}]{ReducibilityLand}, the Spanier--Whitehead dual of $X$ is given by the Thom spectrum of its Spivak normal fibration---precisely the spectrum $\MSF(X)$ as in \cref{def:thom-SF}---thus, there is a canonical equivalence $\bD(\Sigma^\infty X_+)\simeq \MSF(X)$, and the Poincar\'e duality isomorphism may be factored as
\[
    \PD\colon\cR^n(X)\xlongrightarrow{\text{Thom}}\cR^{n-d}(\MSF(X))\xlongrightarrow{\SW}\cR_{d-n}(X).
\]
\cref{thm:A} implies that the right-hand map induces an isomorphism between the cohomological and the homological Atiyah--Hirzebruch spectral sequences. It remains to show that the Thom isomorphism induces an isomorphism between the cohomological Atiyah--Hirzebruch spectral sequences. For the unstable Spivak normal fibration $\SF_\ell(X)\colon X\ra \BG(\ell-d)$, we proceed by identifying the above exact couples defining the spectral sequences via a couple map $\smash{(\Psi, \Phi)\colon (\sA_X^*,\sD_X^*)\ra (\sA_{\SF_\ell(X)}^{*+\ell-d},\sD_{\SF_\ell(X)}^{*+\ell-d})}$ inducing an isomorphism between the spectral sequences. Using \cref{lem:skeleta-orientations}, we recall the following equivalences of spectra. A chosen $\cR$-orientation of the spherical fibration $\SF_\ell(X)\colon X\ra \BG(\ell-d)$ induces an $\cR$-orientation on the restriction to subskeleta. Thus, there is a Thom equivalence
\[
    \TD_c\colon\F(\Sigma^\infty X^{(p)}_+,\cR)\lra \F(\Sigma^\infty\Th(\SF_\ell(X))^{(p+\ell-d)}[d-\ell],\cR)
\]
between mapping spectra. Furthermore, taking cofibres of inclusions of subskeleta, we obtain an induced equivalence 
\[
    \TD_c^\cofib\colon\F(\Sigma^\infty(X^{(p)}/X^{(p-1)}),\cR)\lra \F(\Sigma^\infty(\Th(\SF_\ell(X))^{(p+\ell-d)}/\Th(\SF_\ell(X))^{(p-1+\ell-d)})[d-\ell],\cR)
\]
as follows: By replacing $\F(-,\S)$ with $\F(-,\cR)$ and repeating the proof of \cref{lem:cofib} verbatim, the skeletal filtrations of $X$ and $\Th(\SF_\ell(X))$, as in the definition of the exact couples above, yield fibration sequences on mapping spectra. In particular, we obtain a commutative diagram
\[\begin{tikzcd}
	{\F(\Sigma^\infty X^{(p-1)}_+,\cR)[-1]} & {\F(\Sigma^\infty\Th(\SF_\ell(X))^{(p-1+\ell-d)}[d-\ell],\cR)[-1]} \\
	{\F(\Sigma^\infty(X^{(p)}/X^{(p-1)}),\cR)} & {\F(\Sigma^\infty(\Th(\SF_\ell(X))^{(p+\ell-d)}/\Th(\SF_\ell(X))^{(p-1+\ell-d)})[d-\ell],\cR)} \\
	{\F(\Sigma^\infty X^{(p)}_+,\cR)} & {\F(\Sigma^\infty\Th(\SF_\ell(X))^{(p+\ell-d)}[d-\ell],\cR)} \\
	{\F(\Sigma^\infty X^{(p-1)}_+,\cR)} & {\F(\Sigma^\infty\Th(\SF_\ell(X))^{(p-1+\ell-d)}[d-\ell],\cR)}
	\arrow["{\TD_c}", from=1-1, to=1-2]
	\arrow["{k^*}", from=1-1, to=2-1]
	\arrow["{\kappa^*}", from=1-2, to=2-2]
	\arrow["{\TD_c^\cofib}", from=2-1, to=2-2]
	\arrow["{j^*}", from=2-1, to=3-1]
	\arrow["{\lambda^*}", from=2-2, to=3-2]
	\arrow["{\TD_c}", from=3-1, to=3-2]
	\arrow["{i^*}", from=3-1, to=4-1]
	\arrow["{\iota^*}", from=3-2, to=4-2]
	\arrow["{\TD_c}", from=4-1, to=4-2]
\end{tikzcd}\]
of fibration sequences. By \cref{lem:skeleta-orientations}, we know that the horizontal maps $\TD_c$ are equivalences, and a short application of the Five Lemma together with Whitehead's theorem yields that $\TD_c^\cofib$ is an equivalence as well.
Note that, by the definition of $\cR$-cohomology groups, $\cR^{p+q-1}(X^{(p)})$ is \emph{equal to} $\smash{\pi_{-p-q+1}(\F(\Sigma^\infty X_+^{(p)},\cR))}$, and similarly for the other mapping spectra in the above diagram. Thus, we may define the couple map $(\Psi\coloneqq (\TD_c)_*,\Phi\coloneqq (\TD_c^\cofib)_*)$ as the maps on homotopy groups, induced by the Thom equivalences on subskeleta and cofibres. The morphisms do indeed define a couple map: Briefly, applying homotopy groups to the diagram above yields a commutative diagram---which is exactly what is required for $(\Psi,\Phi)$ to define a couple map. As in the proof of \cref{thm:A}, we suppress an additional suspension isomorphism implicit in the definition of the exact couples, which is justified by naturality of the suspension isomorphism. Since $\TD_c$ and $\TD_c^\cofib$ are equivalences of mapping spectra, $(\Psi,\Phi)$ defines an \emph{isomorphism} between the exact couples, which leads to an isomorphism between the two associated spectral sequences.
Hence, we have shown that the unstable Thom isomorphism 
\[
    -\cup u(\SF_\ell(X))\colon\cR^n(X)\lra\cR^{n+\ell-d}(\Sigma^\infty \Th(\SF_\ell(X))).
\]
induces an isomorphism between the Atiyah--Hirzebruch spectral sequence computing $\cR^*(X)$ and the one computing $\cR^{*+\ell-d}(\Sigma^\infty \Th(\SF_\ell(X)))$.

We now pass to the stable version $\MSF(X)$ as follows. Firstly, we remark that the shifting operator yields an identification $\smash{\cR^{n+\ell-d}(\Sigma^\infty\Th(\SF_\ell(X)))\simeq \cR^{n-d}(\Sigma^\infty \Th(\SF_\ell(X))[-\ell])}$ which induces an isomorphism between the Atiyah--Hirzebruch spectral sequences. This may be seen by passing to mapping spectra; for any spectrum $Y$, there is an equivalence of mapping spectra $\F(Y,\cR)\simeq \F(Y[k],\cR[k])$, which induces a suspension isomorphism $\cR^0(Y)\simeq \cR^k(Y[k])$. In particular, the suspension isomorphism exists on the level of mapping spectra, and hence induces an isomorphism between the Atiyah--Hirzebruch spectral sequences by naturality.
Since $\ell$ is assumed to be large enough that $\SF_\ell(X)$ represents the Spivak normal fibration, there is an equivalence $\Sigma^\infty \Th(\SF_{N+1}(X))\simeq \Sigma^\infty\Th(\SF_N(X))[1]$ of suspension spectra for all $N\geq \ell$. Furthermore, each map in the diagram
\[
    \Sigma^\infty\Th(\SF_\ell(X))[-\ell]\lra\Sigma^\infty\Th(\SF_{\ell+1}(X))[-(\ell+1)]\lra\dots
\]
contravariantly induces an isomorphism between the Atiyah--Hirzebruch spectral sequences computing $\cR^*(-)$. The limit of the spectral sequences is precisely the one computing
\[
    \lim_{N\geq \ell}\cR^{*}(\Sigma^\infty\Th(\SF_N(X))[-N])\simeq\cR^{*}(\colim_{N\geq \ell}\Sigma^\infty\Th(\SF_N(X))[-N])\simeq\cR^{*}(\MSF(X)).
\]
Note that the limit of the spectral sequences exists since every morphism of spectral sequences is an isomorphism. The unstable Thom isomorphisms induce isomorphisms between the Atiyah--Hirzebruch spectral sequence computing $\cR^*(X)$ and the ones computing $\cR^{*+N-d}(\Sigma^\infty \Th(\SF_N(X)))$ for all $N\geq \ell$. By commutativity of the diagram at the end of the proof of \cref{lem:MSF-orientation}, as well as naturality of the Atiyah--Hirzebruch spectral sequence, these isomorphisms are compatible with the stabilisation isomorphisms $\cR^{*+N-d}(\Sigma^\infty \Th(\SF_N(X)))\ra \cR^{*+N+1-d}(\Sigma^\infty \Th(\SF_{N+1}(X)))$ on the level of spectral sequences. Putting this together, we obtain an isomorphism of spectral sequences from the Atiyah--Hirzebruch spectral sequence computing $\cR^*(X)$ to the spectral sequence obtained as the limit over $N\geq \ell$ of the spectral sequences computing $\cR^{*-d}(\Sigma^\infty \Th(\SF_N(X))[-N])$. Once more, since all involved morphisms are isomorphisms, the latter agrees with the Atiyah--Hirzebruch spectral sequence computing $\smash{\lim_{N\geq \ell}\cR^{*-d}(\Sigma^\infty\Th(\SF_N(X))[-N])\simeq \cR^{*-d}(\MSF(X))}$.
In particular, this isomorphism of spectral sequences is induced by the equivalence of spectra
\[
    \F(\Sigma_+^\infty X,\cR)\lra \lim_{N\geq \ell}{}\F(\Sigma^\infty\Th(\SF_N(X))[-(N-d)],\cR)\simeq \F(\MSF(X)[d],\cR),
\]
which in turn is precisely the stable Thom isomorphism $\cR^{*}(X)\ra \cR^{*-d}(\MSF(X))$ upon taking homotopy groups.
\end{proof}

\begin{Ackn}
    I thank Daniel Kasprowski for several enlightening discussions, as well as Oscar Randal-Williams for helpful comments (\cref{rmk:r-w}) on a previous version of this article, and for encouragement to give more details at the end of \cref{sec:PD-Orientation}. Furthermore, I thank the anonymous referee for carefully reading the article and providing many helpful suggestions. I was supported by the University of Southampton, and the UKRI (Grant EP/W524621/1).
\end{Ackn}

\bibliographystyle{amsalpha}
\bibliography{References}

\end{document}